\documentclass[9pt]{article}
\usepackage{latexsym}
\usepackage{xcolor}
\usepackage{amssymb}
\usepackage{MnSymbol}
\usepackage{float}
\restylefloat{table}
\usepackage{booktabs}
\usepackage{mathrsfs}
\usepackage{amsmath}   
\usepackage{graphicx}
\usepackage{geometry}
\usepackage[font=footnotesize]{caption}

\begin{document}
\begin{center}

{\large \textbf{BASKETS AND FIBRED LINKS REALIZING $A_{n}$}}\\
\vspace{20px}
 \small{ LUCAS FERNANDEZ VILANOVA}\footnote{The author is supported by the Swiss National Science Foundation (project no. 178756).}
 \end{center}

\renewcommand{\abstractname}{}
\begin{abstract}
\footnotesize {\normalsize A}BSTRACT. We prove that basket links, whose symmetrized Seifert form is congruent to the Cartan matrix of the simply laced Dynkin diagram $A_{n}$, are isotopic to the torus link $T(2,n+1)$. In addition, we provide examples of links, constructed by plumbing $n$ positive Hopf bands the core curves of which intersect at most once, with symmetrized Seifert form congruent to the Cartan matrix $A_{n}$, that are not isotopic to $T(2,n+1)$.  \par
\end{abstract}

\vspace{10px}
\begin{center}
 1. {\large I}NTRODUCTION\\
\end{center}

\begin{figure}[H]
\includegraphics[scale=0.5]{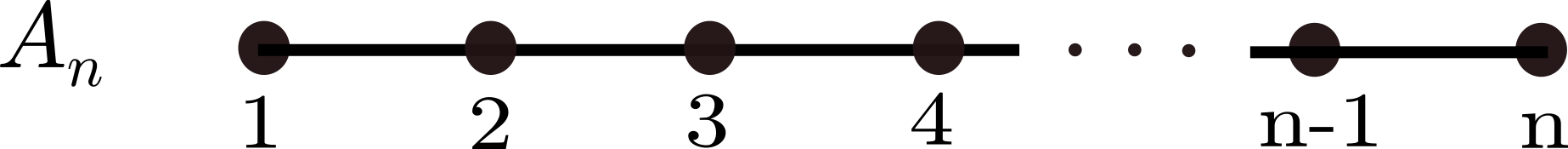}\ \ \ \ \ \ \ \ \ \ \ \ 
\includegraphics[scale=0.17]{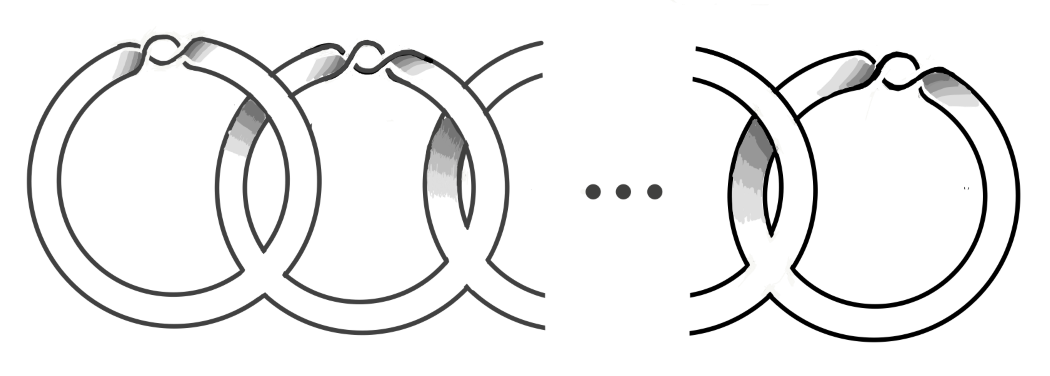}
\centering
\caption{Left: the $A_{n}$ diagram. Right: positive Hopf bands plumbed according to the diagram $A_{n}$.}
 \label{fig:an}
\end{figure}

The torus links $T(2,n+1)$ can be constructed by plumbing $n$ positive Hopf bands according to the diagram $A_{n}$ (see Figure \ref{fig:an}); such construction is an example of what is known as a \textit{positive arborescent Hopf plumbing}. We will denote by $C_{A_{n}}$ the symmetrized Seifert matrix of the link $T(2,n+1)$ (tridiagonal matrix with 2's in the main diagonal and $-1$'s in the upper and lower diagonals). While the congruence class of $C_{A_{n}}$ has been studied by graph theorist in the context of the adjacency matrix, see e.g. \cite{I}, the links realizing such matrices are far from being understood. Initial works in this direction show that these are, for some types of links, precisely the two strand torus links $T(2,n+1)$. It is the case of positive braids as proved by Baader, \cite{S}, and for certain basket links as pointed out by Boileau, Boyer and Gordon, \cite{BBG}. Here, we show that this also works for basket links: \\

\textbf{Theorem 1.1.} \textit{A basket link with $n$ positive Hopf bands and symmetrized Seifert form congruent to $C_{A_{n}}$ is isotopic to a two strand torus link.}\\

\begin{figure}[H]
\includegraphics[scale=0.14]{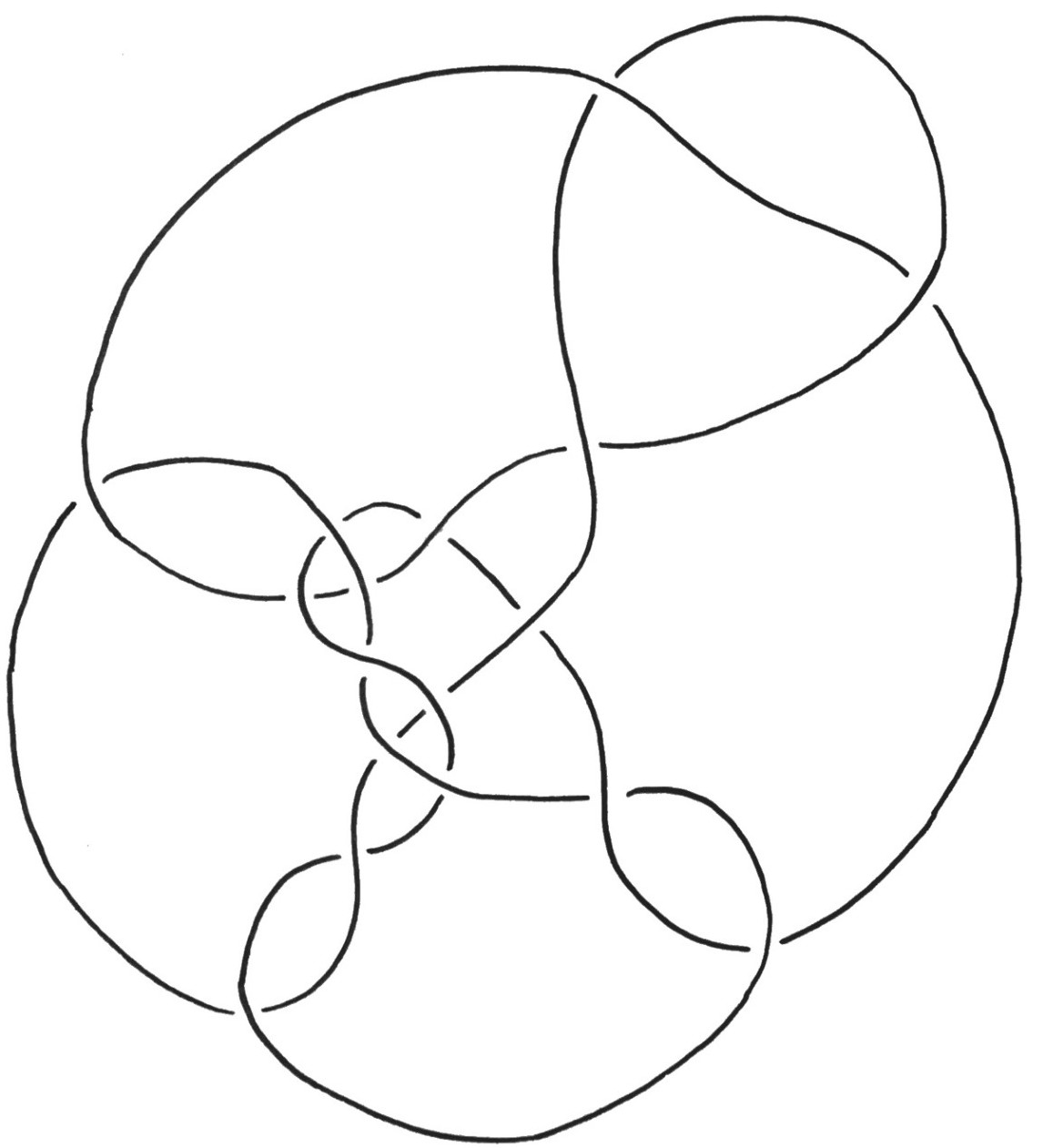}
\centering
\caption{A $17$-crossings fibred knot resulting from plumbing $6$ positive Hopf bands. It has the same Seifert form as $T(2,7)$ yet they are not isotopic. Image plotted with \textit{knotscape}.}
\label{fig:plumb5}
\end{figure}

 F. Misev found that there exists an infinite family of distinct fibred knots having the same Seifert form as the torus knots $T(2,2g+1)$ for any given genus $g
\geq 2$, \cite{FM}. These knots were constructed by plumbing positive Hopf bands such that their core curves intersect more than once. Since for basket links, the core curves of the Hopf bands intersect at most once, this fact together with Theorem 1.1 motivates the question whether we can extend the result of Theorem 1.1. to the more general class of links arising by plumbing positive Hopf bands in which the core curves of the bands intersect at most once. The last section is dedicated to show, via an example, that this is not true. In the same line as in Misev's article, we show how to construct links and knots with the same Seifert form as the Torus links $T(2,n+1)$ but distinct form it, see Figure \ref{fig:plumb5} for an anticipated example. 

These links are interesting given the recent conjecture proposed in \cite{BBG}, stating that fibred strongly quasipositive links such that its cyclic branched cover is an L-space are \textit{simply laced arborescent}, e.g., of the type $T(2,n+1)$. Secondly, these examples also become interesting in the context of the prevailing slice-ribbon conjecture, since Baker proved that if $K$ and $K'$ are two distinct fibred strongly quasipositive knots, then $K\#-K'$ is not ribbon, \cite{KB}. Although we show that the two component links ($n$ odd) that we produce are not smoothly concordant to $T(2,n+1)$, the question remains open for $n$ even.\\

\textbf{1.2. Outline}\\
In the second section we give a short introduction to basket links. In the third section we show that symmetrized Seifert matrices that are congruent to $C_{A_{n}}$ carry a signed graph (a graph with positive or negative edges) with a specific structure; they are what we call \textit{complete-tree graphs}. This uses previous results on spectral graph theory and constitutes the first part of the proof of Theorem 1.1. In the fourth section, we use this to show by sliding bands, that a basket with such intersection graph is isotopic to a two strand torus link, completing the proof. In the last section we explain how to construct links with the same Seifert form as the two strand torus links but non-isotopic to them. The construction method is to plumb an additional band to a given basket link.  \\

\textbf{Acknowledgements}\\
I thank S. Baader for guidance and for proofreading this paper. I also thank the referee for useful comments and suggestions.

\begin{center}
 2. {\large B}ASKETS\\
\end{center}

\begin{figure}[H]
\includegraphics[scale=0.3]{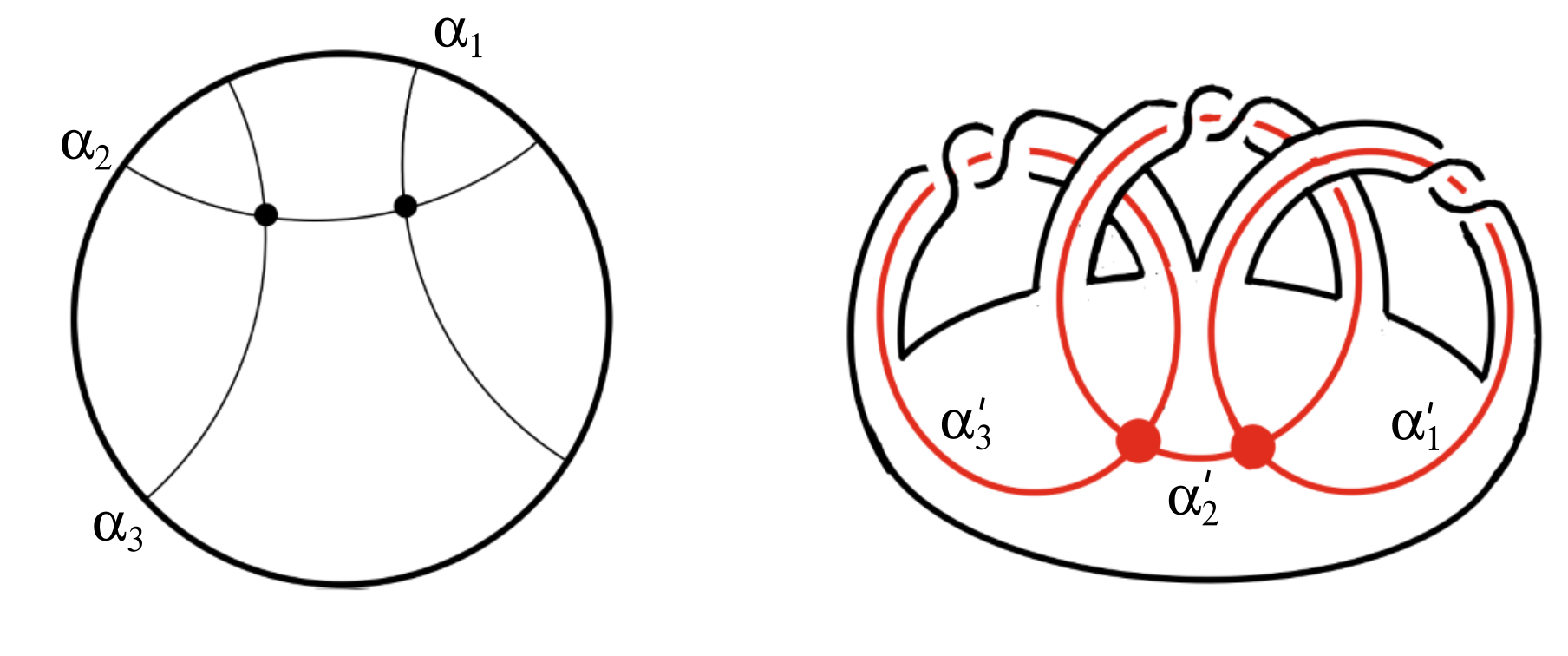}
\centering
\caption{ A chord diagram, left, and its basket surface, right.}
\label{fig:basket}
\end{figure}

Basket links were first defined by Rudolph, \cite{LR}. They are strongly quasipositive fibred links, that include positive braid links. They are constructed by plumbing Hopf bands along the neighborhoods of properly embedded arcs in a disk. The easiest and more intuitive representation of a basket as it appears in \cite{HI}, is by using the so called chord diagrams; a disk in $\mathbb{R}^{2}$ with ordered chords, see for instance Figure \ref{fig:basket}, left. From there, we can build the basket surface by plumbing Hopf bands in the given order and following the convention of plumbing the $i^{th}$ band on bottom of the previous bands, see Figure \ref{fig:basket}, right. 

The incidence graph, $\Gamma$, of a chord diagram has one vertex for each chord and such that two vertices are connected whenever the corresponding chords intersect, e.g, Figure \ref{fig:basket} (left) is a chord diagram with incidence graph $A_{3}$. It is worth mentioning that, looking at Figure \ref{fig:basket} (right), if we change a bottom plumbing by a top plumbing, the basket surface does not change, \cite{LR}. There are incidence graphs, e.g, $A_{n}$, where the order in which we plumb is irrelevant. In fact, in \cite{HI}, they show that this is also true for tree graphs. 

A chord diagram with arcs $\alpha_{1},\dots, \alpha_{n}$ together with an order of plumbing determines uniquely a signed graph, $\Gamma_s$, i.e, a graph in which every edge is assigned a value 1 or -1. We obtain $\Gamma_{s}$  as follows: identically as we did with the incidence graph we associate one vertex for each arc and one edge whenever two arcs intersect. In order to find the sign, let $\alpha_{i}'$ be the core curve of the positive Hopf band plumbed along $\alpha_{i}$, then the sign of the edge $[\alpha_{i},\alpha_{j}]$ is given by the sign of $lk(\alpha_{i}',\alpha_{j}'^{+})+lk(\alpha_{i}'^{+},\alpha_{j}')$, where $\alpha_{i}'^{+}$ is the curve resulting from pushing $a_{i}'$ off in the positive normal direction of the surface. Clearly, the incidence graph is the underlying graph of $\Gamma_{s}$. Moreover, in this basis, the basket surface has symmetrized Seifert matrix $M(\Gamma_s)=2I+A(\Gamma_s)$. Here, $A(\Gamma_s)$ is the adjacency matrix of $\Gamma_{s}$ (where $A(\Gamma_{s})_{ij}=\pm 1$ if the vertices $i$ and $j$ are connected by an edge). Recall from the previous section that 
\[
C_{A_{n}}=\left( \begin{array}{ccccccc}
2&-1&& &  &\\
-1&2&-1\\
 &-1&2&-1\\
 & & &\ddots\\
 & & &-1&2&-1\\
& & & &-1&2\\
\end{array}\right).
\]

The matrix $C_{A_{n}}$ is positive definite, so if $M(\Gamma_s)\cong C_{A_{n}}$ then so is $M(\Gamma_s)$ and the eigenvalues of $A(\Gamma_{s})$ are $>-2$. This motivates the next section where we investigate the structure of such graphs. \\

\begin{center}
 3. {\large A} NOTE ON $A_{n}$ GRAPHS\\
\end{center}

Graphs whose eigenvalues are larger than $-2$ have been studied  and classified by Cameron et al. in connection with the root systems $A_{n}$, $D_{n}$, $E_{6}$, $E_{7}$ and $E_{8}$, \cite{CGSS}. A characterization of these graphs was given by Doob and Cvetkovic, \cite{MD}. Later, in \cite{GJ}, Greaves et al. found a similar characterization for signed graphs. In fact, Ishihara shows that the signed graphs of type $A_{n}$ come with a specific structure, named \textit{Fushimi trees}, \cite{I}. Here, we give an insight in such theories and show an alternative proof to the one of Ishihara.\\

Let $G$ be a finite, connected, simple and signed graph. Denote by $\lambda (G)$ the least eigenvalue of $A(G)$, and by $M(G)$ the matrix $2I+A(G)$. If $\lambda (G)>-2$ we will say that $G$ is definite. 

The root system $A_{n}$ is the set of vectors in $\mathbb{R}^{n+1}$ of the form $\pm (e_{i}-e_{j})$ for $1\leq i< j\leq n+1$. A graph $G$ (signed or not) is said to be represented by a root system if $M(G)=KK^{T}$, where all the rows of $K$ are vectors in the root system, \cite{LW}. The line graph of a graph $G$, denoted as $L(G)$, is a graph in which every vertex represents an edge of $G$, and two vertices are adjacent if and only if their corresponding edges share a vertex. For example, the complete graph $K_{n}$ is the line graph of a star graph with $n+1$ vertices. In \cite{LW}, Theorem 6.3, they show that a graph is represented by the root system $A_{n}$ if and only if it is the line graph of a bipartite graph. It is not hard to show, that if a graph is signed and represented by $A_{n}$, then its underlying graph is the line graph of a bipartite graph. Denote by $|G|$ the underlying graph of a signed graph $G$.\\

\textbf{Lemma 3.1.}\textit{ Let $G$ be a signed graph such that $M(G)$ is congruent to $C_{A_{n}}$, then $|G|$ is the line graph of a tree.}\\

 \textit{Proof:} If $M(G)$ is congruent to $C_{A_{n}}$, then $G$ is represented by the root system $A_{n}$, meaning that $|G|=L(S)$ for a bipartite graph $S$. Now, since the line graph of an $n$-cycle is again a $n$-cycle, if $S$ has an even cycle, then so does $|G|$, but then $M(G)$ is not congruent to $C_{A_{n}}$ by \cite{L}, independently of the signs. Therefore, $S$ must be a tree. \hfill$\square$\\ 

For simplicity, we will call a signed and definite graph of type $L(T)$, for some tree $T$ by \textit{complete-tree graphs} (or \textit{Fushimi trees}). Note that, switching the signs of the edges incident to a given vertex of a graph, $G$, does not change the congruence class of $M(G)$. Two signed graphs are  \textit{switching equivalent} if we can transform one into another by switching signs without changing the congruence class. It is worth mentioning, that a complete-tree graph is switching equivalent to a complete-tree graph with only positive edges, see \cite{I}. \\

\begin{center}
 4. {\large P}ROOF OF THEOREM 1.1.\\
\end{center}

As proved by \cite{BBG} Theorem 9.11, definite graphs of type $A_{n}$ and $K_{m}$ have a unique realization as baskets. Indeed, they show how these two baskets (when $n=m$) are isotopic by sliding bands. Complete-tree graphs have a unique realization as baskets as well (they are a tree-like amalgamation of the former). In this section we will prove by induction on the number of complete subgraphs, that baskets realized by a complete-tree graph are isotopic to one with incidence graph $A_{n}$. Here, we will adopt the notation in \cite{BBG}, and we will say that a leg of size $n$ is an $A_{n}$ graph attached at one of its leaves to a vertex of a complete graph, and that a basket is definite when its symmetrized Seifert form is positive definite.\\

\textbf{Lemma 4.1:}\textit{ Let $B_{s}$ be a definite basket with $s$ bands, the incidence graph of which is $K_{n}$ with $m\leq n$ legs and $n>2$, then $B_{s}$ is isotopic to a basket with incidence graph of type $A_{s}$.}\\

\begin{figure}[H]
\includegraphics[scale=0.45]{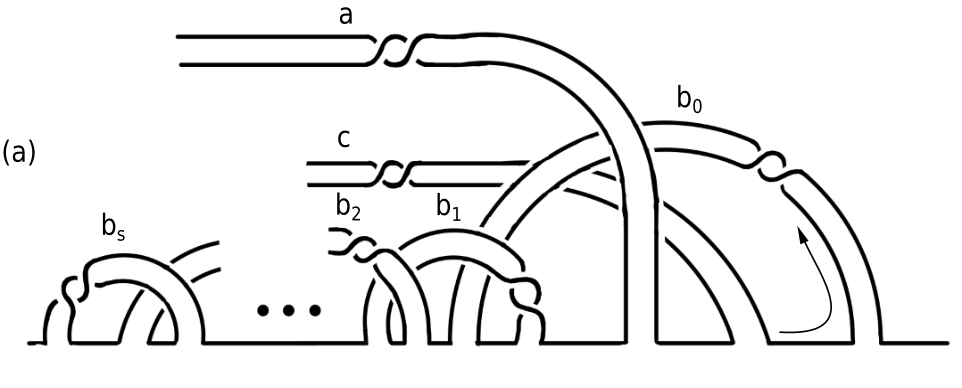}
\includegraphics[scale=0.46]{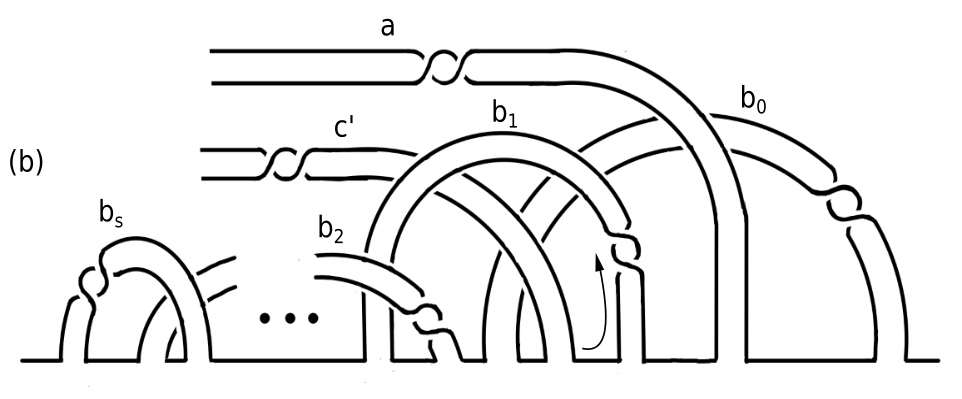}
\includegraphics[scale=0.45]{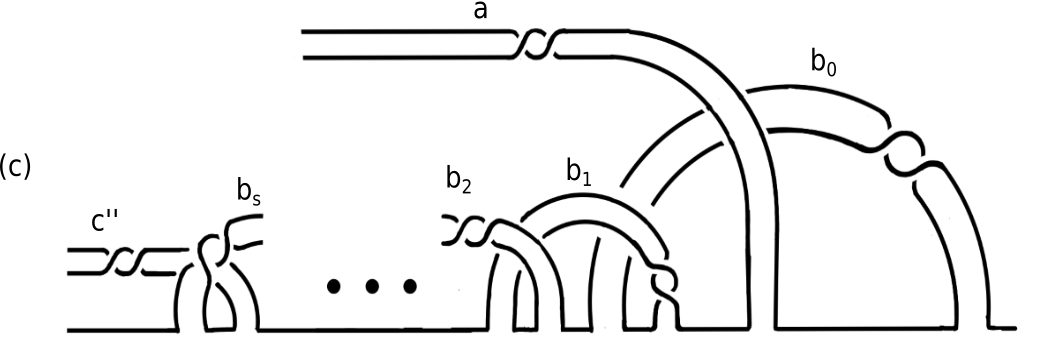}
\centering
\caption{Slide $c$ over $b_{0}$ to obtain $c'$ and repeat.}
 \label{fig:acb}
\end{figure}

\textit{Proof:} First, consider the configuration of positive Hopf bands plumbed in the disk $D$ as it appears in Figure \ref{fig:acb} (a). In terms of the incidence graph the bands $b_{1},\dots, b_{s}$ form a leg of size $s$. Slide $c$ over $b_{0}$ (obtaining the bands in (b)), then $c'$ over $b_{1}$ and continue in this fashion till sliding it over $b_{s-1}$, resulting in (c). 

\begin{figure}[H]
\includegraphics[scale=0.29]{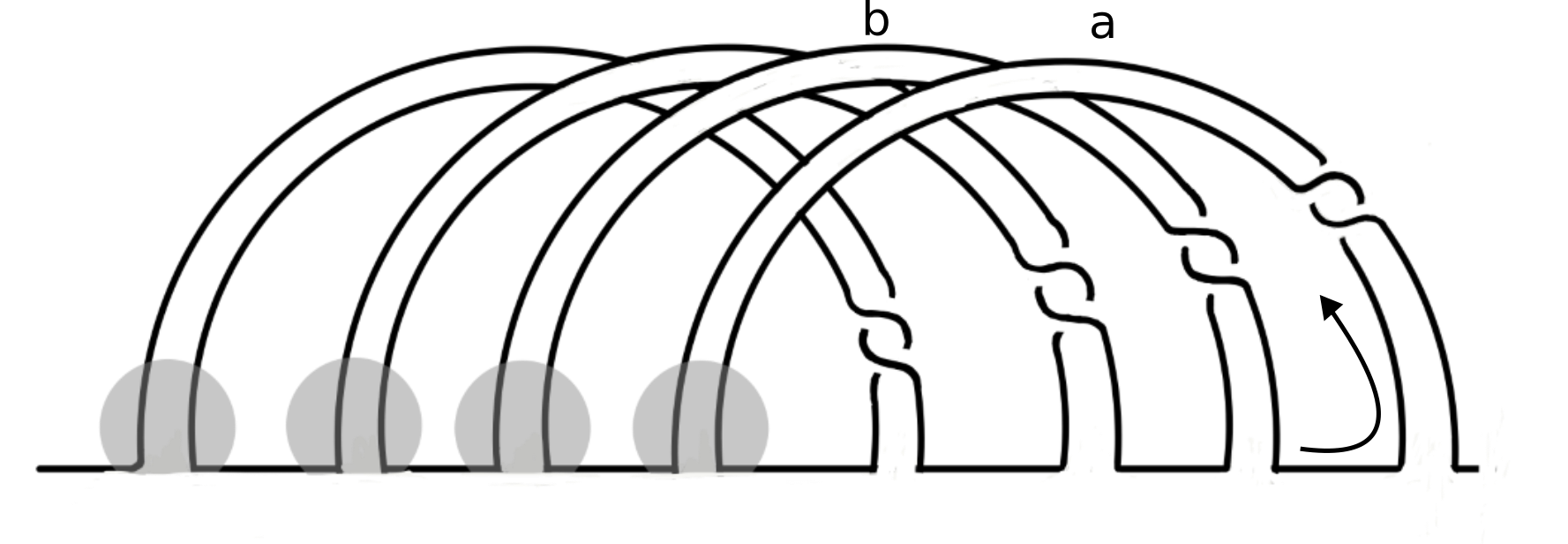}
\centering
\caption{The graph $K_{4}$ with four legs (represented by the shadow areas). By sliding $b$ over $a$ and using the moves in Figure \ref{fig:acb} we obtain $K_{3}$ with three legs.}
 \label{fig:k4}
\end{figure}

Now, consider as an example the complete graph $K_{4}$ with $4$ legs (the cases with less number of legs work similarly) shown in Figure \ref{fig:k4}, where the shadow areas represent where the legs are glued (for simplicity, the legs are not shown in the figure). Here, we assume that the legs extend to the left, the proof works similarly otherwise. Note that sliding $b$ over $a$ and then sliding it over each of the subsequent  bands forming the leg connected to $a$ (using the moves in Figure \ref{fig:acb}) results in $K_{3}$ with $3$ legs.\hfill$\square$\\

\textit{Proof of Theorem 1.1:} If a basket link with $n$ positive Hopf bands has a symmetrized Seifert form congruent to $C_{A_{n}}$, by Lemma 3.1 we know that the incidence graph of such basket is the underlying graph of a complete-tree graph, say $G$. We will proceed by induction on the number of complete graphs in $G$. Observe that, in $G$ there is always at least one outermost complete graph with legs, such that the subgraph of $G$ that results after deleting it, is connected. Since, by Lemma 4.1, the basket with incidence graph the complete graph with legs is isotopic to one with an $A_{n}$ incidence graph, the induction follows.\hfill$\square$\\

\begin{center}
 5. {\large P}LUMBING ON A BASKET\\
\end{center}

In this section we will explain how to construct links, by plumbing positive Hopf bands with symmetrized Seifert form congruent to $C_{A_{n}}$ for $n\geq 5$ which are not the torus link $T(2,n+1)$, and such that the core curves of these Hopf bands intersect at most once. Recall that this differs from the knots studied by F. Misev, \cite{FM}, where core curves intersect more than once. We consider separately the cases in which $n$ is odd and even. We will see that the former case are links composed of two components: one of them is the pretzel link $P(n-3,1,1,1)$ while the other is unknotted. 

\begin{figure}[H]
\includegraphics[scale=0.42]{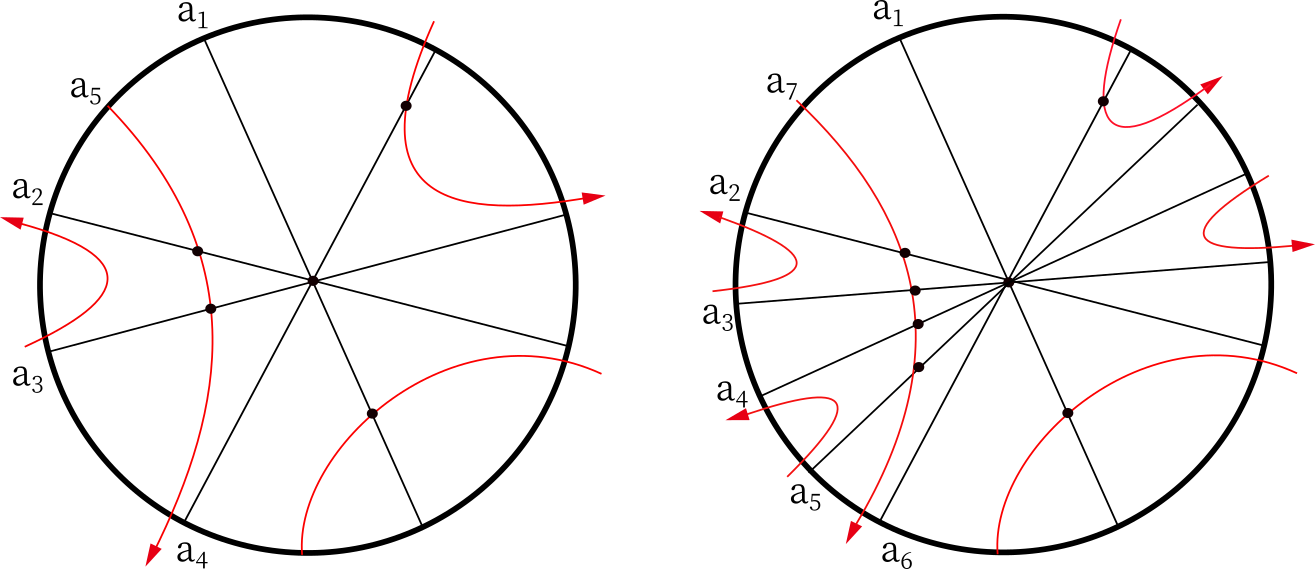}
\centering
\caption{Left: Chord diagram with four arcs $a_{1},a_{2},a_{3},a_{4}$ together with a red arc, $a_{5}$, that travels along the Hopf bands indicated with an arrow. Right: an example with seven Hopf bands.}
\label{fig:plumb1}
\end{figure}

\begin{figure}[H]
\includegraphics[scale=0.23]{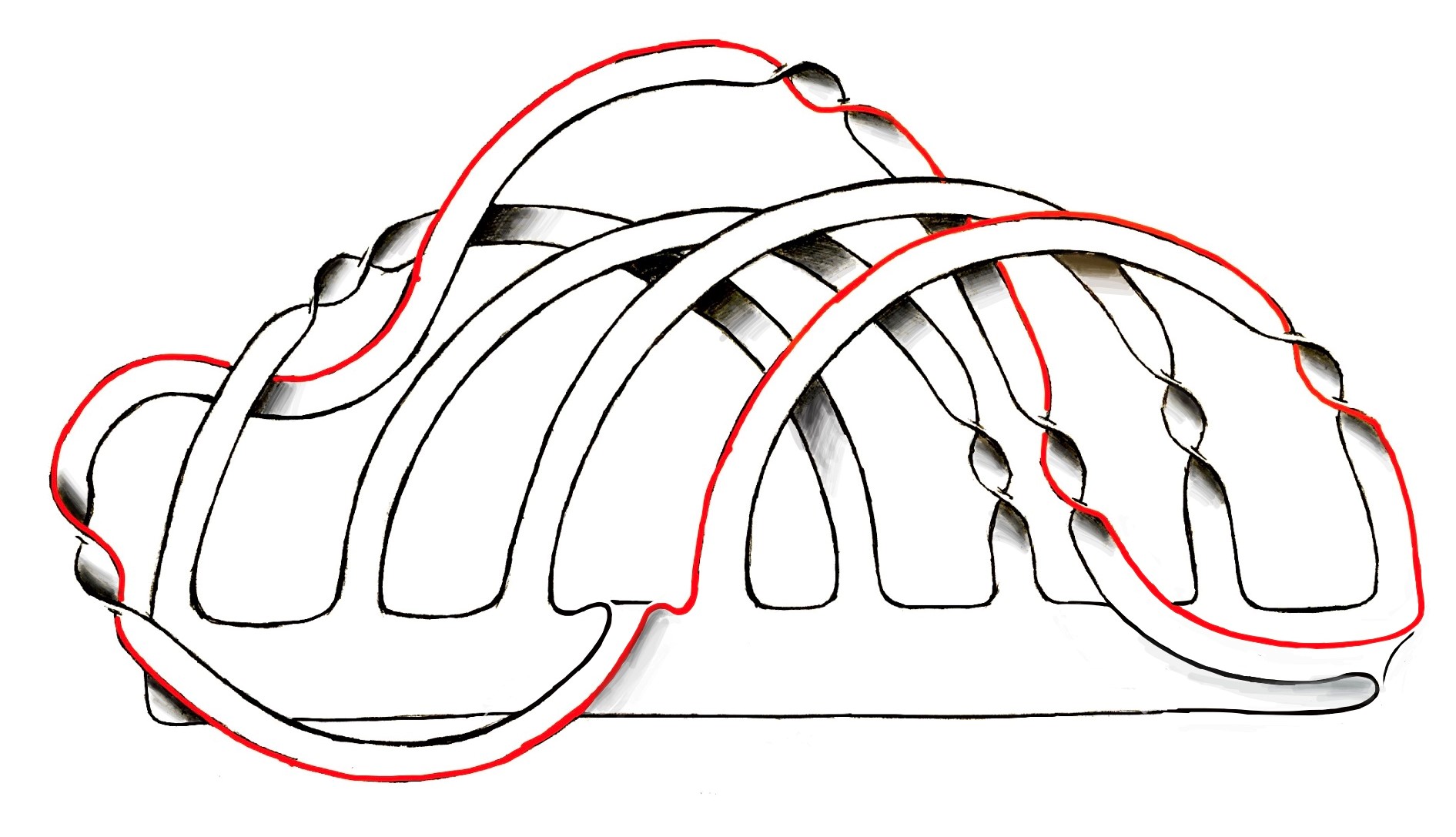}
\centering
\caption{Five positive Hopf bands plumbed according to the diagram in Figure 6, left.}
\label{fig:plumb2}
\end{figure}

 Let us consider the case in which the  number of arcs is odd, say $n=2m+1$ for $m>2$. Consider a disk $D\subset \mathbb{R}^{2}$ with properly embedded arcs $a_{1},\dots,a_{2m}$ (a chord diagram) such that its incidence graph is $K_{2m}$. Let $S_{2m}$ be the surface after plumbing positive Hopf bands on bottom of each other along a neighborhood of the arcs $a_{i}$ for $1<i\leq2m$. Let $\beta$  be a properly embedded arc parallel to $a_{1}$ and intersecting all but $a_{1}$ (note that there are two equivalent possibilities). Then plumb along the arc $a_{2m+1}=T_{a_{2m}}\cdots T_{a_{2}}(\beta)$ (or $T_{a_{2}}\cdots T_{a_{2m}}(\beta)$), where $T_{a}(b)$ stands for the Dehn twist of the curve $a$ along $b$. It is clear that now the core curves of the Hopf bands intersect according to the complete graph $K_{2m+1}$. In Figure \ref{fig:plumb1} there are two examples of this construction: for five and seven Hopf bands. We can check that the Seifert matrix has congruence type $A_{n}$. In fact, this is easier to check on Figure \ref{fig:plumb2}.

\begin{figure}[H]
\includegraphics[scale=0.7]{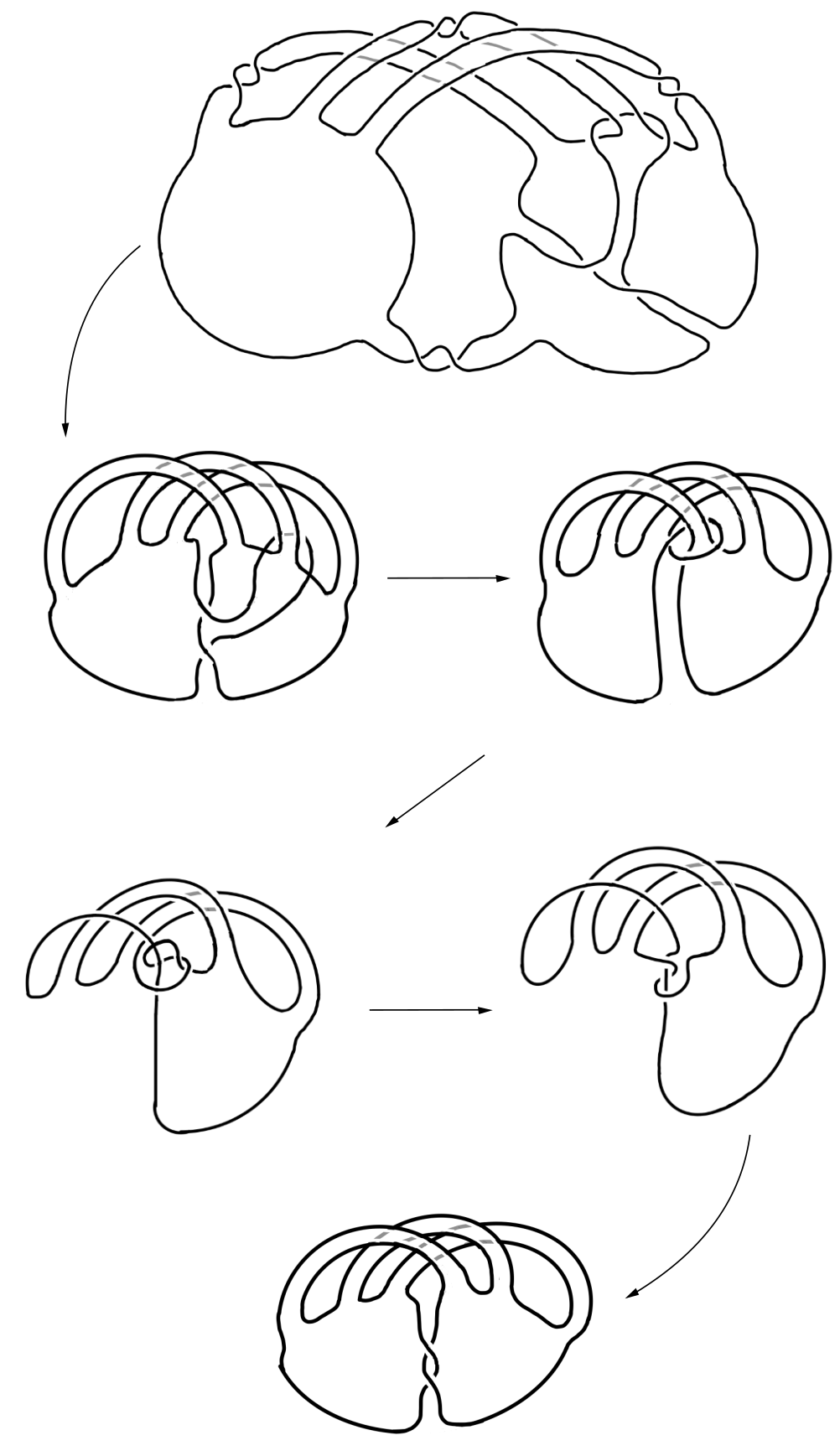}
\centering
\caption{ The component marked with black of the link in Figure 7, and how to transform it into a recognisable $5_{2}$ knot.}
\label{fig:plumb3}
\end{figure}

We will move our attention to the case with five bands, the resulting surface of which, after isotopy, appears in Figure \ref{fig:plumb2}. This is a two component link. The linking number is $3$ and can be easily computed using Figure \ref{fig:plumb2}, where one component is marked with red (unknotted!) and the other with black. The latter is the knot $5_{2}$. Therefore, showing that this cannot be the Torus link $T(2,6)$. In order to see this, if we delete the red component, we are left with a rather complicated knot that we have rearranged (using Reidemeister moves) for our convenience as it appears in Figure \ref{fig:plumb3}, from which is relatively easy to find moves that brings the knot into a a recognisable pretzel knot $P(2,1,1,1)$ (also known as $5_{2}$ or $P(3,1,1)$). In particular, this shows that this link is not smoothly concordant to $T(2,6)$, since $5_{2}$ is not concordant to the unknot. \textit{Snappy}, \cite{SNAPP}, has shown that this is an hyperbolic link, therefore it has a pseudo-Asonov monodromy. 

\begin{figure}[H]
\includegraphics[scale=0.28]{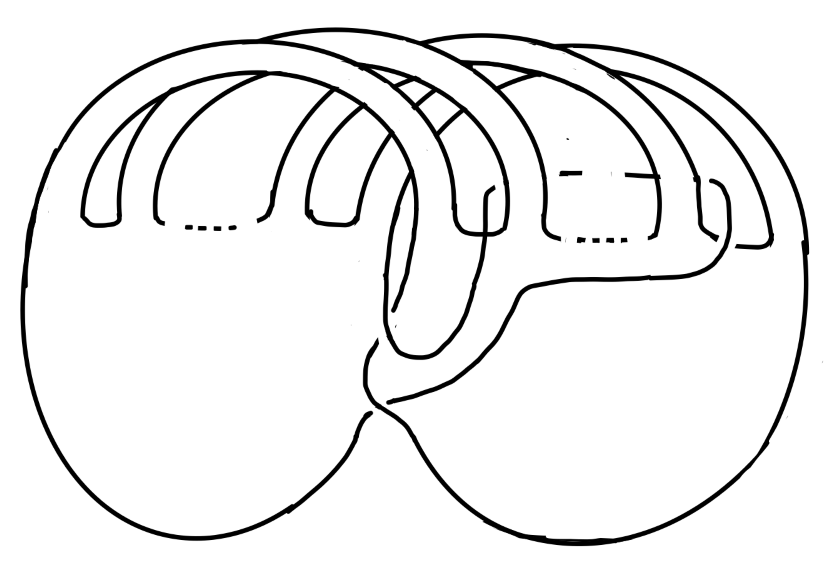}
\centering
\caption{The knotted component for $n>5$.}
\label{fig:plumb4}
\end{figure}

Similarly as we did for the case with five bands, if $n>5$, we still obtain an unknotted component and the other one can be put as in Figure \ref{fig:plumb4}, from there it is not too hard to obtain the pretzel knots $P(n-3,1,1,1)$ in a similar fashion as in Figure \ref{fig:plumb3}.

The case with six Hopf bands can be obtained using a chord diagram with five arcs and incidence graph $K_{5}$. Let $\beta$ be an arc intersecting $a_{2}$, $a_{3}$ and $a_{4}$ (there are two equivalent possibilities), and plumb the sixth band along $a_{6}=T_{a_{4}}T_{a_{3}}T_{a_{2}}(\beta)$ (or $a_{6}=T_{a_{2}}T_{a_{3}}T_{a_{4}}(\beta)$). The resulting knot has $17$ crossings (computed with \textit{Knotscape}, \cite{KS}, see Figure \ref{fig:plumb5}) and its Jones polynomial is different from the one of $T(2,7)$.\\

\textbf{Remark 5.1.} The number of links obtained by positive Hopf plumbing where the core curves intersect at most once times is finite. Indeed, by induction, assume that there are only a finite number of surfaces that can be obtained by plumbing $n$ positive Hopf bands. Plumb another positive Hopf band to that surface. This must be done along a properly embedded arc in this surface. Now, a result from Przytycki shows that the number of essential simple arcs on a punctured surface with Euler characteristic $\chi<0$, that are pairwise non-homotopic and intersect at most once is exactly $2|\chi|(|\chi|+1)$, \cite{PP}. The induction follows.

Mathematics Institute, University of Bern, Sidlerstrasse 5, 3012 Bern, Switzerland\\
\textit{E-mail address:} lucas.fernandez@math.unibe.ch

\end{document}